\documentclass[11pt,reqno]{amsart}
\usepackage{amssymb,latexsym}
\usepackage{graphicx}

\theoremstyle{plain}
\numberwithin{equation}{section}
\newtheorem{thm}{Theorem}[section]
\newtheorem{theorem}[thm]{Theorem}

\newtheorem{proposition}[thm]{Proposition}

\begin{document}

\setcounter{page}{1}

\title{ Number of Compositions  and Convolved Fibonacci numbers }
\author{Milan Janji\'c}
\address{Department of Mathematics and Informatics\\
                Banja Luka University\\
                78000 Banja Luka, Republic of  Srpska\\
                Bosnia and Herzegovina}
\email{agnus@blic.net}
\begin{abstract}
We consider two type of upper Hessenberg matrices which determinants are Fibonacci numbers. Calculating sums of principal minors of the fixed order of the first type  leads us to convolved Fibonacci numbers. Some identities for these  and for Fibonacci numbers are proved.

We also show that numbers of compositions of a natural number with fixed number of ones  appear as coefficients of characteristic polynomial of  a Hessenberg matrix which determinant is  a Fibonacci number. We derive the explicit formula for the number of such compositions, in terms of convolutions of  Fibonacci numbers.
\end{abstract}

\maketitle

\section{Introduction}
Different numbers may be represented as determinants of Hessenberg matrices. There are several such matrices which determinants are Fibonacci numbers. We consider two
types of such matrices. It occurs that coefficients of the characteristics polynomial of matrices of the first type, that is, the sums of its principal minors, are convolved Fibonacci numbers in  the sense of \cite{riordan}. From this fact we derive several identities. It is also shown that convolved Fibonacci numbers are coefficients of expansion of Fibonacci polynomial $f_{n+1}(x-1)$ in terms of powers of $x.$

Coefficients of characteristics polynomial of matrices of the seconde type produce numbers $c(n,k)$ of compositions of the natural number $n$ with exactly $k$ ones. We derive the explicit formula for $c(n,k)$ in terms of convolutions of Fibonacci numbers.

The following known result about upper Hessenberg
matrices will be used in the paper.
\begin{theorem}\label{tt} Let $a_1,p_{i,j},\;( i\leq j)$ be arbitrary elements of  a commutative ring $R,$ and let the sequence $a_1,a_2,\ldots$ be defined by:
\begin{equation}\label{niz}a_{n+1}=\sum_{i=1}^np_{i,n}a_i,\;(n=1,2,\ldots).\end{equation}

If
\begin{equation}A_n=
\left[\begin{array}{llllll}
p_{1,1}&p_{1,2}&p_{1,3}&\cdots&p_{1,n-1}&p_{1,n}\\
-1&p_{2,2}&p_{2,3}&\cdots&p_{2,n-1}&p_{2,n}\\
0&-1&p_{3,3}&\cdots&p_{3,n-1}&p_{3,n}\\
\vdots&\vdots&\vdots&\ddots&\vdots&\vdots\\
0&0&0&\cdots&p_{n-1,n-1}&p_{n-1,n}\\
0&0&0&\cdots&-1&p_{n,n}\end{array}\right]\end{equation} then
\[a_{n+1}=a_1\det A_n,\;(n=1,2,\ldots).\]
\end{theorem}

As usual, $f_{-1},f_0,f_1,\ldots$ are Fibonacci numbers. Note that $f_{-1}=1,f_0=0.$
\section{Convolved Fibonacci Numbers}
It seems that convolved Fibonacci numbers are first appeared in the classical Riordan's book \cite{riordan}. Here we shall obtain these numbers in another way.
Consider the following upper Hessenberg matrix:
\[F_n=\left[\begin{array}{rrrrrr}
1&1&0&\cdots&0&0\\
-1&1&1&\cdots&0&0\\
0&-1&1&\cdots&0&0\\
\vdots&\vdots&\vdots&\ddots&\vdots&\vdots\\
0&0&0&\cdots&1&1\\
0&0&0&\cdots&-1&1\end{array}\right].\]

From Theorem \ref{tt} we obtain the following well-known result: \[\det F_n=f_{n+1},\;(n=1,2,\ldots).\]

Principal minors of $F_n$ are some convolutions of Fibonacci numbers. For example, the minor obtained by deleting the $i$th row and
column of $F_n$ is obviously  equal to $f_{i}\cdot f_{n-i+1}.$ It follows
that the principal minor $M(i_1,i_2,\ldots,i_k)$ of $F_n$ obtained by
deleting rows and columns with indices $1\leq
i_1<i_2<\cdots<i_k\leq n$ is
\[M(i_1,i_2,\ldots,i_k)=f_{i_1}\cdot f_{i_2-i_1}\cdots
f_{i_k-i_{k-1}}\cdot f_{n-i_k+1}.\]
Thus, the following proposition holds:
\begin{proposition} The sums $S_{n-k},\;(k=0,1,2,\ldots,n-1)$ of principal minors of the order $n-k$ of $F_n$ is
\[S_{n-k}=\sum_{j_1+j_2+\cdots+j_{k+1}=n-k}f_{j_1+1}f_{j_2+1}\cdots
f_{j_{k+1}+1},\] where the sum is taken over  $j_t\geq
0,\;(t=1,2,\ldots,k+1).$ \end{proposition}
In the sense of \cite{riordan} the sum on the right side of the preceding equation is called convolved Fibonacci number and is denoted by $f^{(k+1)}_{n-k+1}.$ Hence,

\[S_{n-k}=f^{(k+1)}_{n-k+1}.\]

The characteristic matrix of $F_n$ has the form:
\[\left[\begin{array}{rrrrrr}
x-1&1&0&\cdots&0&0\\
-1&x-1&1&\cdots&0&0\\
0&-1&1&\cdots&0&0\\
\vdots&\vdots&\vdots&\ddots&\vdots&\vdots\\
0&0&0&\cdots&x-1&1\\
0&0&0&\cdots&-1&x-1\end{array}\right].\]
It follows that for characteristic polynomial $p_{n}(x)$ of $F_n$ holds $p_n(x)=f_{n+1}(x-1),$ where $f_{n+1}(x)$ are Fibonacci polynomials.

Therefore, we have the following:
\begin{proposition}
Let $f_{n}(x)$ be the Fibonacci polynomial. Then \[q_{n+1}(x-1)=\sum_{k=0}^n(-1)^{n-k}f^{(k+1)}_{n-k+1}x^k.\]
\end{proposition}

From the well-known explicit formula for Fibonacci polynomials we obtain the following:
\begin{proposition} Let $f^{(k+1)}_{n-k+1}$ be convolved Fibonacci numbers. Then
\[f^{(k+1)}_{n-k+1}=\sum_{i=0}^{\lfloor\frac{n-k}{2}\rfloor}{n-i\choose i}{n-2i\choose k}.\]
\end{proposition}

Since $q_{n+1}(1)=f_{n+1}$ we obtain the following identity for Fibonacci numbers:
\begin{proposition} For Fibonacci numbers $f_n,\;(n=1,\ldots)$ we have
\[f_{n+1}=(-1)^n\sum_{k=0}^n\sum_{i=0}^{\lfloor\frac{n-k}{2}\rfloor}(-2)^k{n-i\choose i}{n-2i\choose k}.\]
\end{proposition}

The matrix $F_n$ is obviously invertible. It is easy to find its inverse.
Denote by $M(i,j)$ the cofactor  of the element from the $i$th row and the $j$th column of $F_n.$ It follows that:
\[M(i,j)=f_i\cdot f_{n-j+1},\;(i\leq j);\;M(i,j)=(-1)^{i+j}f_j\cdot f_{n-i+1},\;(i>j).\]
From a known property of  associated matrices we obtain an expression of powers of Fibonacci numbers in terms of some of its convolutions.
\begin{proposition}
If $M_n=(M_{ij})_{n\times n}$ is the associated matrix of $F_n$ then \[\det M=f_{n+1}^{n-1}.\]
\end{proposition}

\section{Number of Compositions of a Natural Number with Fixed numbers of Ones}
We shall now consider another type of  Hessengerg matrices which determinants are Fibonacci numbers. This will lead us to the number of compositions of a natural number with fixed number of ones.
\begin{proposition}\label{pp2}
Let $G_n$ be the matrix of order $n$ defined by:
\[G_n=\left[\begin{array}{rrrrrr}
0&1&1&\cdots&1&1\\
-1&0&1&\cdots& 1&1\\
0&-1&0&\cdots&1&1\\
\vdots&\vdots&\vdots&\ddots&\vdots&\vdots\\
0&0&0&\cdots&0&1\\
0&0&0&\cdots&-1&0\end{array}\right].\]
Then \[\det(G_n)=f_{n-1}.\]
\end{proposition}
\begin{proof}
The equation $\det(G_1)=0,\;\det G_2=1$ are obvious. According to (\ref{tt}) we have $\det(G_{n+1})=\det(G_0)+\det(G_1)+\cdots+\det(G_{n-1}).$
This recurrence shows that the proposition is true.

\end{proof}
\begin{proposition}\label{pp1}
Let $S_{n-k},\;(0\leq k\leq n)$ be the sum of all minors of order $n-k$ of $G_n.$  Then
\[S_{n-k}=\sum_{j_1+j_2+\cdots+j_{k+1}=n-2k+1}f_{j_1}f_{j_2}\cdots f_{j_{k+1}},\]
where the sum is taken over  $j_t\geq -1,\;(t=1,2,\ldots,k+1).$
\end{proposition}
\begin{proof}

Let $M(i_1,i_2,\ldots,i_k)$ be the minor of order $n-k$ obtained by
deleting rows and columns with indices $1\leq
i_1<i_2<\cdots<i_k\leq n.$ Then
\[M=\det(G_{i_1-1})\cdot\det(G_{i_2-i_1-1})\cdots\det(G_{i_k-i_{k-1}-1})\cdot\det(G_{n-i_k}).\]
Applying Proposition \ref{pp2} we obtain
\[M=f_{i_1-2}f_{i_2-i_1-2}\cdots f_{i_k-i_{k-1}-2} f_{n-i_k-1}.\]
It follows that
\[S_{n-k}=\sum_{1\leq
i_1<i_2<\cdots<i_k\leq n}f_{i_1-2}f_{i_2-i_1-2}\cdots f_{i_k-i_{k-1}-2} f_{n-i_k-1},\] that is,
\[S_{n-k}=\sum_{j_1+j_2+\cdots+j_{k+1}=n-2k-1}f_{j_1}f_{j_2}\cdots f_{j_{k+1}},\]
where the sum is taken over all  $j_t\geq -1,\;(t=1,2,\ldots,k+1).$
\end{proof}
We shall finished with an explicit formula for the number  $c(n,k),$ of compositions of $n$ with exactly $k$ ones.
\begin{theorem} Let $n$ be a positive integer.  Then for $0\leq k\leq n$ we have 
\begin{equation}\label{cnk}c(n,k)=\sum_{j_1+j_2+\cdots+j_{k+1}=n-2k+1}f_{j_1}f_{j_2}\cdots f_{j_{k+1}},\end{equation}
where the sum is taken over $j_t\geq -1,\;(t=1,2,\ldots,k+1).$
\end{theorem}
\begin{proof}
We use the induction with respect of $k.$ 
For $k=0,$ according to Proposition \ref{pp2} and  Proposition \ref{pp1}, we have $c(n,0)=S_n=f_{n-1}.$
Hence, the theorem is true for $k=0.$

Assume that theorem is true for $k-1.$ The greatest  value of $j_{k+1}$ is $n-k-1,$ and is obtained
for $j_1=j_2=\cdots=j_k=-1.$
We thus may write $(\ref{cnk})$ in the form:
\[c(n,k)=\sum_{j=-1}^{n-k-1}f_j\sum_{j_1+j_2+\cdots+j_{k}=n-2k+1-j}f_{j_1}f_{j_2}\cdots f_{j_{k}}.\]
By the induction hypothesis we have
\begin{equation}\label{cnk1}c(n,k)=\sum_{j=-1}^{n-k-1}f_j\cdot c(n-j-2,k-1).\end{equation}
Let $(i_1,i_2,\ldots)$ be a composition of $n$ with exactly $k$ ones,  and let the first $1$ occurs at  the $p$th place.
Then $(i_1,i_2,\ldots,i_{p-1})$ is a composition of $j+1$ with no ones, and $(i_{p+1},\ldots)$ is a composition of $n-j-2$ with exactly $k-1$ ones, 
where $j+1=i_1+\cdots+i_{p-1}.$ 
It follows that the number of compositions in which the first $1$ is on the $p$th place is  $f_j\cdot c(n-j-2,k-1).$ This is a term in the sum on the right side of (\ref{cnk1}).

For $p=1$ we have $j=-1$ which produces the first term in (\ref{cnk1}). If ones  are on the last $k$ places of a composition  then
$j+1=n-k$ that gives $f_{n-k-1}c(k-1,k-1)=f_{n-k-1}$ which is the last term in the sum of (\ref{cnk1}).
\end{proof}

Note that triangle of $c(n,k),(n=0,1,\ldots;\;k=0,1,\ldots,n)$ appears as A105422   in Sloane's Encyclopedia \cite{sloane}. Also, $c(n,k)$ appear in another context in the paper \cite{papa}, where they are introduced trough a generating function and a recurrence relation. It is also shown that that $c(n,k)$
is the number of $n$ length bit strings beginning with $0,$ having $k$ singles.

\medskip

\noindent AMS Classification Numbers: 11B39, 15A36

\end{document}